\documentclass[final]{siamltex}
\usepackage{amsfonts}
\usepackage{amssymb,pstcol,pstricks}

\topskip  \input{tcilatex}

\begin{document}

\title{Grid polygons from permutations and their enumeration by the kernel
method}
\author{Toufik Mansour\thanks{%
Department of Mathematics, University of Haifa, Haifa 31905, Israel (\texttt{%
toufik@math.haifa.ac.il}).} \and Simone Severini\thanks{%
Department of Mathematics and Department of Computer Science, University of
York, York YO10 5DD, U.K (\texttt{ss54@york.ac.uk}).}}
\maketitle

\begin{abstract}
A \emph{grid polygon} is a polygon whose vertices are points of a grid. We
define an injective map between permutations of length $n$ and a subset of
grid polygons on $n$ vertices, which we call \emph{consecutive-minima
polygons}. By the \emph{kernel method}, we enumerate sets of permutations
whose consecutive-minima polygons satisfy specific geometric conditions. We
deal with $2$-variate and $3$-variate generating functions involving
derivatives, cases which are not routinely solved by the kernel method.
\end{abstract}



\pagestyle{myheadings} \thispagestyle{plain} 
\markboth{T. MANSOUR
AND S. SEVERINI}{Grid polygons from permutations}

\section{Introduction and examples}

The enumeration of permutations that satisfy certain constraints has
recently attracted interest (see \emph{e.g.} \cite%
{AA,Mi03,GM,M1,M2,Va,We95,Z,ZZ}). In particular, permutation patterns have
been extensively studied over the last decade (see for instance \cite{GM}
and references therein). The tools involved in these works include
generating trees (with either one or two labels), combinatorial approaches,
recurrences relations, enumeration schemes, scanning elements algorithms, 
\emph{etc}.

Permutations are traditionally associated to a number of combinatorial and
algebraic objects, like matrices, trees, posets, graphs, \emph{etc}. (see 
\emph{e.g.} \cite{go,st}). The purpose of this work is to begin a study of
the interplay between permutations and polygons. A practical motivation
comes from computational geometry, where the complexity of algorithms for
polygons is an important subject \cite{co}. Of course, it is intuitive that
imposing combinatorial constraints on geometric objects commonly reduces
generality; on the other side, techniques from combinatorics may provide a
fertile background for the design of algorithms, even if the analysis is
restricted to toy-cases.

Here we associate permutations to a subset of the grid polygons and
enumerate sets of permutations whose polygons satisfy specific geometric
conditions. Clearly, there are many potential ways to associate permutations
to polygons. Each ways presumably has a special feature which helps to
underline some particular property of the permutations. If we want to keep a
one-to-one correspondence, this arbitrariness is materialized in two points:

\begin{itemize}
\item Of all possible polygons associated to a given permutation, we choose
the one with a fixed extremal property, for example, the polygon with
minimum area or perimeter.

\item We decide how to construct a polygon according to some chosen rule.
The rule should guarantee the association of each permutation to a single
polygon, unequivocally.
\end{itemize}

\begin{figure}[h]
\begin{center}
\begin{pspicture}(.5,.5)(4.5,3.5)
\psgrid[unit=0.5,subgriddiv=1,griddots=10,gridlabels=6pt](1,1)(7,7)
\pscircle*(0.5,2){0.05}\pscircle*(1,2.5){0.05}\pscircle*(1.5,1){0.05}
\pscircle*(2,1.5){0.05}\pscircle*(2.5,.5){0.05}\pscircle*(3,3.5){0.05}\pscircle*(3.5,3){0.05}
\end{pspicture}
\begin{pspicture}(.5,.5)(4,3.5)
\psgrid[unit=0.5,subgriddiv=1,griddots=10,gridlabels=6pt](1,1)(7,7)
\pscircle*(0.5,2){0.05}\pscircle*(1,2.5){0.05}\pscircle*(1.5,1){0.05}
\pscircle*(2,1.5){0.05}\pscircle*(2.5,.5){0.05}\pscircle*(3,3.5){0.05}\pscircle*(3.5,3){0.05}
\psline(0.5,2)(1.5,1)(2.5,.5)(3.5,3)(3,3.5)(2,1.5)(1,2.5)(.5,2)
\end{pspicture}
\end{center}
\caption{The grid $L_{4523176}$ and the consecutive-minima polygon $%
P_{4523176}$.}
\label{flat}
\end{figure}
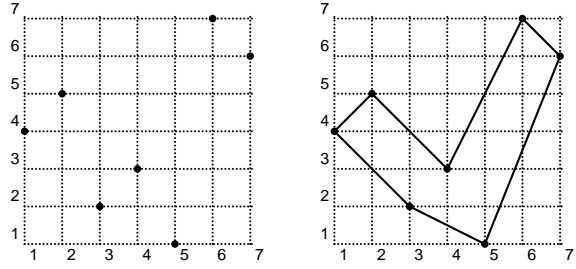

We opt here for the second approach, as it is formalized in what follows. A 
\emph{grid of side }$n$ is an $n\times n$ array containing $n^{2}$ points, $%
n $ in each row and each column. The distance between two closest points in
the same row or in the same column is usually taken to be $1$ unit. A \emph{%
permutation of length }$n$ is a complete ordering of the elements of the set 
$[n]=\{1,...,n\}$. We associate a grid of side $n$, denoted by $L_{\pi }$,
to a permutation $\pi $ of length $n$. If the permutation takes $i$ to $%
j=\pi _{i}$, we mark\ the point $(i,j)$ of the grid, that is the point in
the row $i$ and the column $j$. For example, the grid $L_{\pi }$ represented
in Figure~\ref{flat} is associated to the permutation $\pi =4523176$.

A \emph{grid polygon} on $n$ vertices is a polygon whose vertices are $n$
points of a grid. A \emph{permutation polygon} on $n$ vertices is a grid
polygon with the following properties: the side of the grid is $n$; in every
row and every column of the grid there is one and only one vertex of the
polygon. It is intuitive to observe that a permutation can be associated to
more than one polygon depending on how we connect the marked points of the
grid. We need to fix some terminology. Let $L_{\pi }$ be a grid of side $n$
of a permutation $\pi $.

\begin{itemize}
\item A point $(i,j)$ is said to be a \emph{left-right minimum} of $L_{\pi }$
if there is no point $(i^{\prime },j^{\prime })$ of $L_{\pi }$ such that $%
i^{\prime }<i$ and $j^{\prime }<j$.

\item A point $(i,j)$ is said to be a \emph{right-left minimum} of $L_{\pi }$
if there is no point $(i^{\prime },j^{\prime })$ of $L_{\pi }$ such that $%
i^{\prime }>i$ and $j^{\prime }<j$.

\item A point $(i,j)$ is said to be a \emph{source} of $L_{\pi}$ if either $%
i=1$, $i=n$, or $(i,j)$ is not a left-right minimum or a right-left-minimum.
\end{itemize}

We say that two points $(i,j)$ and $(i^{\prime },j^{\prime })$ of the grid $%
L_{\pi }$ (resp. of the set of left-right minima, right-left minima,
sources) are \emph{consecutive} if there is no point $(a,b)$ in the set of
(resp. left-right minima's, right-left minima's, sources) such that $%
i<a<i^{\prime }$ or $i^{\prime }<a<i$. For example, the left-right-minima of 
$L_{4523176}$ are $(1,4)$, $(3,2)$ and $(5,1)$; the right-left-minima are $%
(7,6)$ and $(5,1)$; the sources are $(1,4)$, $(2,5)$, $(4,3)$, $(6,7)$, and $%
(7,6)$ (see Figure~\ref{flat}). 
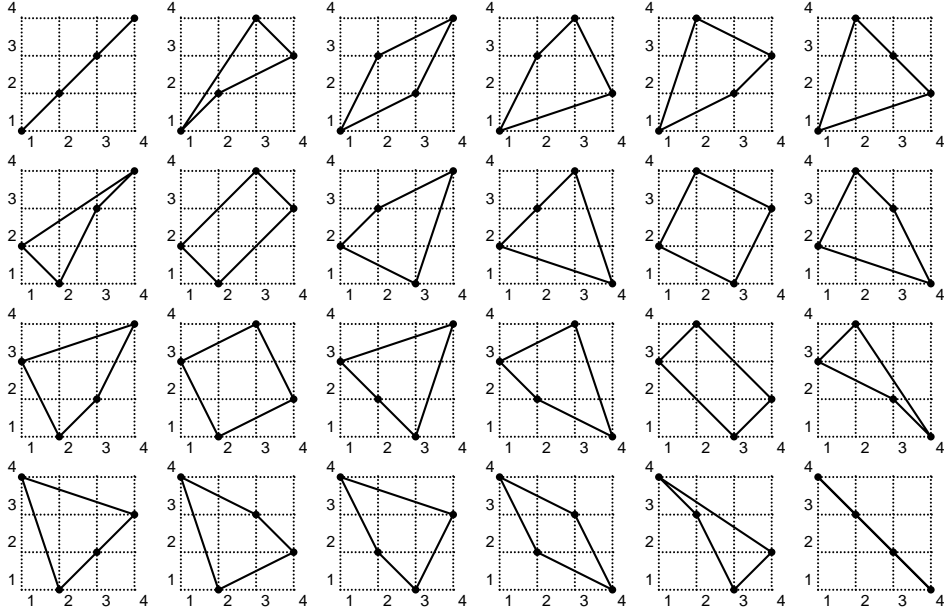
\begin{figure}[h]
\begin{center}
\begin{pspicture}(.5,.5)(2.5,2.5)
\psgrid[unit=0.5,subgriddiv=1,griddots=10,gridlabels=6pt](1,1)(4,4)
\pscircle*(0.5,0.5){0.05}\pscircle*(1,1){0.05}\pscircle*(1.5,1.5){0.05}\pscircle*(2,2){0.05}
\psline(.5,.5)(2,2)
\end{pspicture}
\begin{pspicture}(.5,.5)(2.5,2.5)
\psgrid[unit=0.5,subgriddiv=1,griddots=10,gridlabels=6pt](1,1)(4,4)
\pscircle*(0.5,0.5){0.05}\pscircle*(1,1){0.05}\pscircle*(1.5,2){0.05}\pscircle*(2,1.5){0.05}
\psline(.5,.5)(1,1)(2,1.5)(1.5,2)(.5,.5)
\end{pspicture}
\begin{pspicture}(.5,.5)(2.5,2.5)
\psgrid[unit=0.5,subgriddiv=1,griddots=10,gridlabels=6pt](1,1)(4,4)
\pscircle*(0.5,0.5){0.05}\pscircle*(1,1.5){0.05}\pscircle*(1.5,1){0.05}\pscircle*(2,2){0.05}
\psline(.5,.5)(1.5,1)(2,2)(1,1.5)(.5,.5)
\end{pspicture}
\begin{pspicture}(.5,.5)(2.5,2.5)
\psgrid[unit=0.5,subgriddiv=1,griddots=10,gridlabels=6pt](1,1)(4,4)
\pscircle*(0.5,0.5){0.05}\pscircle*(1,1.5){0.05}\pscircle*(1.5,2){0.05}\pscircle*(2,1){0.05}
\psline(.5,.5)(2,1)(1.5,2)(1,1.5)(.5,.5)
\end{pspicture}
\begin{pspicture}(.5,.5)(2.5,2.5)
\psgrid[unit=0.5,subgriddiv=1,griddots=10,gridlabels=6pt](1,1)(4,4)
\pscircle*(0.5,0.5){0.05}\pscircle*(1,2){0.05}\pscircle*(1.5,1){0.05}\pscircle*(2,1.5){0.05}
\psline(.5,.5)(1.5,1)(2,1.5)(1,2)(.5,.5)
\end{pspicture}
\begin{pspicture}(.5,.5)(2.5,2.5)
\psgrid[unit=0.5,subgriddiv=1,griddots=10,gridlabels=6pt](1,1)(4,4)
\pscircle*(0.5,0.5){0.05}\pscircle*(1,2){0.05}\pscircle*(1.5,1.5){0.05}\pscircle*(2,1){0.05}
\psline(.5,.5)(2,1)(1.5,1.5)(1,2)(.5,.5)
\end{pspicture}
\par
\begin{pspicture}(.5,.5)(2.5,2.5)
\psgrid[unit=0.5,subgriddiv=1,griddots=10,gridlabels=6pt](1,1)(4,4)
\pscircle*(0.5,1){0.05}\pscircle*(1,0.5){0.05}\pscircle*(1.5,1.5){0.05}\pscircle*(2,2){0.05}
\psline(.5,1)(1,.5)(1.5,1.5)(2,2)(.5,1)
\end{pspicture}
\begin{pspicture}(.5,.5)(2.5,2.5)
\psgrid[unit=0.5,subgriddiv=1,griddots=10,gridlabels=6pt](1,1)(4,4)
\pscircle*(0.5,1){0.05}\pscircle*(1,.5){0.05}\pscircle*(1.5,2){0.05}\pscircle*(2,1.5){0.05}
\psline(.5,1)(1,.5)(2,1.5)(1.5,2)(.5,1)
\end{pspicture}
\begin{pspicture}(.5,.5)(2.5,2.5)
\psgrid[unit=0.5,subgriddiv=1,griddots=10,gridlabels=6pt](1,1)(4,4)
\pscircle*(0.5,1){0.05}\pscircle*(1,1.5){0.05}\pscircle*(1.5,.5){0.05}\pscircle*(2,2){0.05}
\psline(.5,1)(1.5,.5)(2,2)(1,1.5)(.5,1)
\end{pspicture}
\begin{pspicture}(.5,.5)(2.5,2.5)
\psgrid[unit=0.5,subgriddiv=1,griddots=10,gridlabels=6pt](1,1)(4,4)
\pscircle*(0.5,1){0.05}\pscircle*(1,1.5){0.05}\pscircle*(1.5,2){0.05}\pscircle*(2,.5){0.05}
\psline(.5,1)(2,.5)(1.5,2)(1,1.5)(.5,1)
\end{pspicture}
\begin{pspicture}(.5,.5)(2.5,2.5)
\psgrid[unit=0.5,subgriddiv=1,griddots=10,gridlabels=6pt](1,1)(4,4)
\pscircle*(0.5,1){0.05}\pscircle*(1,2){0.05}\pscircle*(1.5,.5){0.05}\pscircle*(2,1.5){0.05}
\psline(.5,1)(1.5,.5)(2,1.5)(1,2)(.5,1)
\end{pspicture}
\begin{pspicture}(.5,.5)(2.5,2.5)
\psgrid[unit=0.5,subgriddiv=1,griddots=10,gridlabels=6pt](1,1)(4,4)
\pscircle*(0.5,1){0.05}\pscircle*(1,2){0.05}\pscircle*(1.5,1.5){0.05}\pscircle*(2,.5){0.05}
\psline(.5,1)(2,.5)(1.5,1.5)(1,2)(.5,1)
\end{pspicture}
\par
\begin{pspicture}(.5,.5)(2.5,2.5)
\psgrid[unit=0.5,subgriddiv=1,griddots=10,gridlabels=6pt](1,1)(4,4)
\pscircle*(0.5,1.5){0.05}\pscircle*(1,0.5){0.05}\pscircle*(1.5,1){0.05}\pscircle*(2,2){0.05}
\psline(.5,1.5)(1,.5)(1.5,1)(2,2)(.5,1.5)
\end{pspicture}
\begin{pspicture}(.5,.5)(2.5,2.5)
\psgrid[unit=0.5,subgriddiv=1,griddots=10,gridlabels=6pt](1,1)(4,4)
\pscircle*(0.5,1.5){0.05}\pscircle*(1,.5){0.05}\pscircle*(1.5,2){0.05}\pscircle*(2,1){0.05}
\psline(.5,1.5)(1,.5)(2,1)(1.5,2)(.5,1.5)
\end{pspicture}
\begin{pspicture}(.5,.5)(2.5,2.5)
\psgrid[unit=0.5,subgriddiv=1,griddots=10,gridlabels=6pt](1,1)(4,4)
\pscircle*(0.5,1.5){0.05}\pscircle*(1,1){0.05}\pscircle*(1.5,.5){0.05}\pscircle*(2,2){0.05}
\psline(.5,1.5)(1,1)(1.5,.5)(2,2)(.5,1.5)
\end{pspicture}
\begin{pspicture}(.5,.5)(2.5,2.5)
\psgrid[unit=0.5,subgriddiv=1,griddots=10,gridlabels=6pt](1,1)(4,4)
\pscircle*(0.5,1.5){0.05}\pscircle*(1,1){0.05}\pscircle*(1.5,2){0.05}\pscircle*(2,.5){0.05}
\psline(.5,1.5)(1,1)(2,.5)(1.5,2)(.5,1.5)
\end{pspicture}
\begin{pspicture}(.5,.5)(2.5,2.5)
\psgrid[unit=0.5,subgriddiv=1,griddots=10,gridlabels=6pt](1,1)(4,4)
\pscircle*(0.5,1.5){0.05}\pscircle*(1,2){0.05}\pscircle*(1.5,.5){0.05}\pscircle*(2,1){0.05}
\psline(.5,1.5)(1.5,.5)(2,1)(1,2)(.5,1.5)
\end{pspicture}
\begin{pspicture}(.5,.5)(2.5,2.5)
\psgrid[unit=0.5,subgriddiv=1,griddots=10,gridlabels=6pt](1,1)(4,4)
\pscircle*(0.5,1.5){0.05}\pscircle*(1,2){0.05}\pscircle*(1.5,1){0.05}\pscircle*(2,.5){0.05}
\psline(.5,1.5)(1.5,1)(2,.5)(1,2)(.5,1.5)
\end{pspicture}
\par
\begin{pspicture}(.5,.5)(2.5,2.5)
\psgrid[unit=0.5,subgriddiv=1,griddots=10,gridlabels=6pt](1,1)(4,4)
\pscircle*(0.5,2){0.05}\pscircle*(1,0.5){0.05}\pscircle*(1.5,1){0.05}\pscircle*(2,1.5){0.05}
\psline(.5,2)(1,.5)(1.5,1)(2,1.5)(.5,2)
\end{pspicture}
\begin{pspicture}(.5,.5)(2.5,2.5)
\psgrid[unit=0.5,subgriddiv=1,griddots=10,gridlabels=6pt](1,1)(4,4)
\pscircle*(0.5,2){0.05}\pscircle*(1,.5){0.05}\pscircle*(1.5,1.5){0.05}\pscircle*(2,1){0.05}
\psline(.5,2)(1,.5)(2,1)(1.5,1.5)(.5,2)
\end{pspicture}
\begin{pspicture}(.5,.5)(2.5,2.5)
\psgrid[unit=0.5,subgriddiv=1,griddots=10,gridlabels=6pt](1,1)(4,4)
\pscircle*(0.5,2){0.05}\pscircle*(1,1){0.05}\pscircle*(1.5,.5){0.05}\pscircle*(2,1.5){0.05}
\psline(.5,2)(1,1)(1.5,.5)(2,1.5)(.5,2)
\end{pspicture}
\begin{pspicture}(.5,.5)(2.5,2.5)
\psgrid[unit=0.5,subgriddiv=1,griddots=10,gridlabels=6pt](1,1)(4,4)
\pscircle*(0.5,2){0.05}\pscircle*(1,1){0.05}\pscircle*(1.5,1.5){0.05}\pscircle*(2,.5){0.05}
\psline(.5,2)(1,1)(2,.5)(1.5,1.5)(.5,2)
\end{pspicture}
\begin{pspicture}(.5,.5)(2.5,2.5)
\psgrid[unit=0.5,subgriddiv=1,griddots=10,gridlabels=6pt](1,1)(4,4)
\pscircle*(0.5,2){0.05}\pscircle*(1,1.5){0.05}\pscircle*(1.5,.5){0.05}\pscircle*(2,1){0.05}
\psline(.5,2)(1,1.5)(1.5,.5)(2,1)(.5,2)
\end{pspicture}
\begin{pspicture}(.5,.5)(2.5,2.5)
\psgrid[unit=0.5,subgriddiv=1,griddots=10,gridlabels=6pt](1,1)(4,4)
\pscircle*(0.5,2){0.05}\pscircle*(1,1.5){0.05}\pscircle*(1.5,1){0.05}\pscircle*(2,.5){0.05}
\psline(.5,2)(1,1.5)(2,.5)(1.5,1)(.5,2)
\end{pspicture}
\end{center}
\caption{The consecutive-minima polygons on $4$ vertices.}
\label{fsq4}
\end{figure}

\begin{definition}
A \emph{consecutive-minima polygon} (in what follows just
\emph{polygon}) of a permutation $\pi$, denoted by $P_{\pi}$, is a
permutation polygon in which two vertices $a=(i,j)$ and
$b=(i^{\prime},j^{\prime})$, $i<i^{\prime}$, are connected if one of
the following conditions is satisfied:

\begin{itemize}
\item $a$ and $b$ are consecutive left-right minima of $L_{\pi}$;

\item $a$ and $b$ are consecutive right-left-minima of $L_{\pi}$;

\item $a$ and $b$ are consecutive sources of $L_{\pi}$.
\end{itemize}

In such a context, $(a,b)$ is called an \emph{edge} of $P_{\pi }$.
\end{definition}

For example, the polygon $P_{\pi }$ for all $\pi \in S_{4}$ is represented
in Figure~\ref{fsq4}. In the next sections, we will deal with several
questions about the number of polygons on $n$ vertices that satisfy a
certain set of conditions. In order to do so, we first need to give some
further definition. Let $P$ be a polygon, an edge $((i,j),(i^{\prime
},j^{\prime }))$, $i<i^{\prime }$, of $P$ is said to be \emph{increasing}
(resp. \emph{decreasing}) if $j<j^{\prime }$ (resp. $j>j^{\prime }$). A 
\emph{path} of $P$ is a sequence $(a_{0},a_{1}),(a_{1},a_{2}),\ldots
,(a_{s-1},a_{s})$ of edges of $P$. A \emph{face} of $P$ is either a maximal
path of increasing edges or a maximal path of decreasing edges. For example,
there are exactly $2$, $2$, $16$, $4$ polygons on $4$ vertices of exactly
one, two, three, and four faces, respectively. This can be observed in
Figure~\ref{fsq4}.

\begin{definition}
A polygon is said to be \emph{$k$-faces} if it has exactly $k$
faces. In particular, in the case $k=3$, the polygon is called
triangular, and in case $k=4$ the polygon is called square.
\end{definition}

We will present an explicit formula for the number of $k$-faces polygons on $%
n$, where $k=2,3,4$. It seems to be a challenging question to find an
explicit formula for any $k$.

The technique considered in this paper makes use of generating functions to
convert recurrence relations to functional equations. These are then solved
by the \emph{kernel method} as described in \cite{B02}. It may be
interesting to remark that the kernel method is a routine approach when
dealing with $2$-variate generating functions. However, for functional
equations with more than two variables there is no systematic approach.
Bousquet-M\'{e}lou \cite{Mi03} enumerates four different pattern avoiding
classes of permutations, by using the kernel method with $3$-variate
generating functions. We suggest here another class (namely, the square
permutations), to which corresponds a functional equation defining $3$%
-variate generating functions. Interstingly, such permutations are not
immediately related to pattern avoidance. Among the other techniques, we
remark the use of a two variable functional equation involving a derivative
(see Theorem 4, below), something that does not appear to be common in
enumerative combinatorics.

The remainder of the paper is composed of five sections. In Section 2, we
make some general observations about consecutive-minima polygons. We
characterize convex polygons and enumerate polygons on $n$ vertices with
maximum number of faces. In Section 3, 4 and 5, we enumerate 2-faces,
3-faces and 4-faces polygons, respectively. Section 6 is a list of open
problems.

\section{Some general observations}

\subsection{Convexity}

A polygon is \emph{convex} if the internal angle formed at each vertex is
smaller than $180{{}^\circ}$. Give a sequence $a_{1},a_{2},\ldots ,a_{n}$,
we say that the subsequence $a_{i_{1}},\ldots ,a_{i_{m}}$ with $%
i_{1}<i_{2}<\cdots <i_{m}$ is \emph{fast-growing} if 
\[
\frac{a_{i_{j}}-a_{i_{j+1}}}{i_{j}-i_{j+1}}<\frac{a_{i_{j+1}}-a_{i_{j+2}}}{%
i_{j+1}-i_{j+2}}, 
\]
for any $j=1,2,\ldots ,m-2$, and \emph{slow-growing} if 
\[
\frac{a_{i_{j}}-a_{i_{j+1}}}{i_{j}-i_{j+1}}>\frac{a_{i_{j+1}}-a_{i_{j+2}}}{%
i_{j+1}-i_{j+2}}, 
\]
for any $j=1,2,\ldots ,m-2$.

The consecutive-minima polygon $P_{\pi}$ is convex if and only if

\begin{itemize}
\item the subsequence of left-right-minima of $\pi $ is fast-growing;

\item the subsequence of right-left-maxima of $\pi $ is fast-growing;

\item the subsequence $L_{\pi }$ of the sources of $\pi $ is slow-growing.
\end{itemize}

\subsection{Number of faces}

The number of $1$-face polygons on $n$ vertices is exactly $2$, that is the
polygons corresponding to the two permutations $12\ldots n$ and $n\ldots 1$.
The number of different shapes of consecutive-minima polygon on $n$ vertices
is exactly $n$. To clarify this observation, let $P$ be any
consecutive-minima polygon on $n$ vertices, such that $k$ is maximal if each
face is a segment connecting two vertices, and thus $k\leq n$. It is not
difficult to show that there exists at least one $k$-face consecutive-minima
polygon for any $k=1,2,\ldots ,n$. Let $\phi ^{2k}=2436587\ldots
(2k-2)(2k-3)(2k)(2k+1)(2k+2)\ldots n1(2k-1)$ and $\phi ^{2k+1}=2436587\ldots
(2k)(2k-1)1(2k+1)(2k+2)\ldots n$ be two permutations of length $n$ for all $%
k\geq 1$, then we can see that $P_{\phi ^{k}}$ has exactly $k$-faces. Hence,
for any $n\geq 1$ we have $n$ different shapes of consecutive-minima
polygons on $n$ vertices.

What can be said about permutations with maximum number of faces?\ Let $\pi $
be a permutation of length $n\geq 3$. Since the maximum number of faces of $%
P_{\pi }$ is $n$, one of the following holds:

\begin{enumerate}
\item $\pi $ is an alternating permutation ($\pi $ is said to be \emph{%
alternating} if either $\pi _{1}>\pi _{2}<\pi _{3}>\pi _{4}<\cdots \pi _{n}$
or $\pi _{1}<\pi _{2}>\pi _{3}<\pi _{4}>\cdots \pi _{n}$) such that $\pi
_{1}=1$ and $\pi _{n}=2$;

\item $\pi $ is an alternating permutation such that $\pi _{1}=2$ and $\pi
_{n}=1$;

\item Removing the letter $\pi _{i}=1$, $2\leq i\leq n-1$, from $\pi $ then 
\[
(\pi _{1}-1)\ldots (\pi _{i-1}-1)(\pi _{i+1}-1)\ldots (\pi _{n}-1)
\]
is a permutation satisfying either (1) or (2).
\end{enumerate}

It is not hard to see that the number of alternating permutations of length $%
n$ (see \cite[A000111]{JS} and references therein) satisfying either (1) or
(2) is exactly $E_{n-2}$ if $n$ is odd, otherwise it is $0$ ($E_{n}$ is the
number of alternating permutations on length $n$). Hence, we can state the
following result.

\begin{proposition}
The number of polygons on $n$ vertices with maximum number of faces
($n$ faces) is given by $2E_{n-2}$ if $n$ odd, and $2(n-2)E_{n-3}$
if $n$ even.
\end{proposition}

\section{Enumeration of two faces polygons}

\label{sec22}

A permutation is said to be \emph{parallel} if its polygon has exactly two
faces. For instance, there are $2$ parallel permutations of length $4$,
namely $1324$ and $4231$. We denote the set of all parallel permutations of
length $n$ by $\mathcal{P}_{n}$. Given $a_{1},a_{2},\ldots ,a_{d}\in \mathbb{%
N}$, we define 
\[
p_{n;a_{1},a_{2},\ldots ,a_{d}}=\#\{\pi _{1}\pi _{2}\ldots \pi _{n}\in 
\mathcal{P}_{n}\mid \pi _{1}\pi _{2}\ldots \pi _{d}=a_{1}a_{2}\ldots
a_{d}\}, 
\]
and we denote the cardinality of the set $\mathcal{P}_{n}$ by $p_{n}$. The
main result of this section can be formulated as follows.

\begin{theorem}
\label{pp1} The number of parallel permutations of length $n$ is $$\frac{2}{%
n-1}{{2n-4}\choose{n-2}}-2,$$ for all $n\geq2$.
\end{theorem}

\begin{proof}
First, let us enumerate the permutations of length $n$ that begin at
letter $1$, having polygon with at most two faces. From the
definitions we have that $$p_{n;1}=p_{n;1,2}+p_{n;1,3}+\cdots
+p_{n;1,n}.$$ Besides, for all $a=3,4,\ldots ,n$,
$$p_{n;1,a}=p_{n;1,a,2}+\sum_{j=a+1}^{n}p_{n;1,a,j}=p_{n-1;1,a-1}
+\sum_{j=a+1}^{n}p_{n-1;1,j-1}=\sum_{j=a-1}^{n-1}p_{n-1;1,j},$$ with
the initial conditions $p_{n;1,2}=p_{n-1;1}$. To solve the
recurrence relation of the sequence $p_{n;1,a}$, we need to define
$p_{n}(v)=\sum_{a=2}^{n}p_{n;1,a}v^{a-2}$. Thus, multiplying the
above recurrence relation by $v^{a-2}$ and summing over
$a=3,4,\ldots ,n$ we obtain that
$$
p_{n}(v)=p_{n-1}(1)+\sum_{a=3}^{n}\sum_{j=a-1}^{n-1}p_{n-1;1,j}v^{a-2},
$$
which is equivalent to
$$p_{n}(v)=p_{n-1}(1)+\frac{v}{1-v}(p_{n-1}(1)-vp_{n-1}(v)),$$
for $n\geq 4$. Let $p(v;x)=\sum_{n\geq 3}p_{n}(v)x^{n}$. Multiplying
the above recurrence relation with $x^{n}$ and summing over all
possible $n\geq 4 $, by using the initial condition $p_{3}(v)=1$, we
obtain the following functional equation:
$$p(v;x)=x^{3}+\frac{x}{1-v}(p(1;x)-v^{2}p(v;x)).$$ This type of
equation can be solved using the kernel method \cite{B02}.
Substitute $v=C(x)=\frac{1-\sqrt{1-4x}}{2x}$ in the above functional
equation to get $p(1;x)=x^{3}C^{2}(x)=x^{2}(C(x)-1)$. So, the number
of permutations of length $n$ that begin at letter $1$ and whose
polygon has at most two faces is exactly $c_{n-2}$, for $n\geq 3$,
where $c_{m}=\frac{1}{m+1}{{2m}\choose{m}}$ is the $m$-th Catalan
number. Hence, from the fact there exists only one permutation of
length $n$ whose polygon has exactly one face, namely $12\ldots n$,
we get that the number of parallel permutations of length $n$ that
begin at letter $1$ is $c_{n-2}-1$ for all $n\geq 2$. By making use
of the fact that each parallel permutation of length $n$ can begin
at either $1$ or $n$, we obtain that the number of parallel
permutations of length $n$ is $2(c_{n-2}-1)$, for all $n\geq 2$.
\end{proof}

We now generalize the above enumeration as follows. A polygon is said to be 
\emph{$m$-isolated polygon} if it is a grid polygon contains exactly $m$
isolated points. A permutation $\pi $ is said to be \emph{$m$-parallel} if
begins at letter $1$ and by deleting $m$ vertices from $L_{\pi }$ gives a
polygon with at most two faces. Table \ref{tab1} includes the number of $m$%
-parallel permutations of length $n$ for $m\leq 5$ and $n\leq 12$. 
\begin{table}[ht]
\begin{center}
{\small 
\begin{tabular}{c|llllllllllllll}
$m/n$ & $0$ & $1$ & $2$ & $3$ & $4$ & $5$ & $6$ & $7$ & $8$ & $9$ & $10$ & $%
11$ & $12$ &  \\ \hline
$0$ & $1$ & $1$ & $1$ & $1$ & $2$ & $5$ & $14$ & $42$ & $132$ & $429$ & $%
1430 $ & $4862$ & $16796$ &  \\ 
$1$ & $0$ & $0$ & $0$ & $0$ & $0$ & $1$ & $8$ & $46$ & $232$ & $1093$ & $%
4944 $ & $21778$ & $94184$ &  \\ 
$2$ & $0$ & $0$ & $0$ & $0$ & $0$ & $0$ & $2$ & $26$ & $220$ & $1527$ & $%
9436 $ & $54004$ & $292704$ &  \\ 
$3$ & $0$ & $0$ & $0$ & $0$ & $0$ & $0$ & $0$ & $6$ & $112$ & $1275$ & $%
11384 $ & $87556$ & $608064$ &  \\ 
$4$ & $0$ & $0$ & $0$ & $0$ & $0$ & $0$ & $0$ & $0$ & $24$ & $596$ & $8638$
& $95126$ & $880828$ &  \\ 
$5$ & $0$ & $0$ & $0$ & $0$ & $0$ & $0$ & $0$ & $0$ & $0$ & $120$ & $3768$ & 
$66938$ & $882648$ & 
\end{tabular}%
}
\end{center}
\caption{Number of $m$-parallel permutations of length $n$.}
\label{tab1}
\end{table}

We denote the set of all $m$-parallel permutations of length $n$ by $%
\mathcal{P}_{n}^{m} $. Given $a_{1},a_{2},\ldots ,a_{d}\in \mathbb{N}$, we
define 
\[
p_{n,m;a_{1},a_{2},\ldots ,a_{d}}=\#\{\pi _{1}\pi _{2}\ldots \pi_{n}\in 
\mathcal{P}_{n}^{m}\mid \pi _{1}\pi _{2}\ldots \pi_{d}=a_{1}a_{2}\ldots
a_{d}\}, 
\]
The cardinality of the set $\mathcal{P}_{n}^{m}$ is denoted by $p_{n,m}$.

\begin{theorem}
\label{ppm1} Let $p(v;q,x)=\sum_{n\geq4}\sum_{a=2}^n\sum_{m=0}^n
p_{n,m;1,a}q^mv^{a-2}x^{n}$ be the ordinary generating function for
the sequence $p_{n,m;1,a}$. Then
$$\left(1+\frac{xv^{2}}{1-v}\right)r(v;q,x)-v^2xq\frac{d}{dv}r(v;q,x) =x^3+%
\frac{x}{1-v}r(1;q,x).$$
\end{theorem}

\begin{proof}
Define $p_{n;a_{1},a_{2},\ldots ,a_{d}}=\sum_{m\geq
0}p_{n,m;a_{1},\ldots ,a_{d}}q^{m}$ for any $n$ and $d$. From the
definitions, we have that
$$p_{n}=p_{n;1}=p_{n;1,2}+p_{n;1,3}+\cdots+p_{n;1,n}.$$
Besides, for all $a=3,4,\ldots ,n-1$,
$$\begin{array}{ll}
p_{n;1,a} & =p_{n;1,a,2}+\sum\limits_{j=a+1}^{n}p_{n;1,a,j}+%
\sum\limits_{j=3}^{a-1}p_{n;1,a,j} \\
& =p_{n-1;1,a-1}+\sum\limits_{j=a+1}^{n}p_{n-1;1,j-1}+q(a-3)p_{n-1;1,a-1}=q(a-3)p_{n-1;1,a-1}+\sum\limits_{j=a-1}^{n-1}p_{n-1;1,j},%
\end{array}$$
with the initial conditions $p_{n;1,2}=p_{n-1;1}$ and $p_{n;1,n}=0$.
To solve the recurrence relation of the sequence $p_{n;1,a}$, we
need to define the generating function
$p_{n}(v;q)=\sum_{a=2}^{n}p_{n;1,a}v^{a-2}$. Thus,
multiplying the above recurrence relation by $v^{a-2}$ and summing over $%
a=3,4,\ldots ,n-1$ we obtain that for all $n\geq 5$,
$$p_{n}(v;q)=p_{n-1}(1)+\frac{v}{1-v}(p_{n-1}(1;q)-vp_{n-1}(v;q))+v^{2}q\frac{d%
}{dv}(p_{n-1}(v;q)),$$ with the initial condition $p_{3}(v;q)=1$.
Let $p(v;q,x)=\sum_{n\geq 3}p_{n}(v;q)x^{n}$. If we multiply the
above recurrence relation by $x^{n}$ and we sum over all $n\geq 5$,
we then obtain the requested functional equation.
\end{proof}

\begin{theorem}
\label{ppm2} Let $p_m(v;x)$ be the coefficient of $q^m$ in the
ordinary
generating function $p(v;q,x)$, that is, $p_m(v;x)=\sum_{n\geq3}%
\sum_{a=2}^np_{n,m;1,a}q^mv^{a-2}x^{n}$ (define $p_{-1}(v;x)=0$).
Then
$$\left(1+\frac{xv^{2}}{1-v}\right)p_m(v;x)=v^2x\frac{d}{dv}p_{m-1}(v;x)
\delta_{m,0}x^3+\frac{x}{1-v}p_{m}(1;x),$$ where $\delta_{m,0}=1$ if
$m=0$ and $\delta_{m,0}=0$ otherwise.
\end{theorem}

Theorem \ref{ppm2} provides an algorithm for finding $p_{m}(v;x)$ for any
given $m\geq 0$, since we consider a functional equation with one variable,
which can be solved using the kernel method. It is important to observe that
we cannot just substitute $v=C(x)=\frac{1-\sqrt{1-4x}}{2x}$ in the
functional equation of Theorem \ref{ppm2}, since the generating function $%
\frac{d}{dv}p_{m-1}(v;x)$ is possibly not defined at $v=C(x)$. We need then
the following result, to solve this kind of recurrence relation for given $m$%
.

\begin{theorem}
\label{ppm3} For any $m\geq 0$, the ordinary generating function
$p_{m}(v;x)$ can be written as $$\frac{p_{m}^{\prime
}(v;x)}{(1-v+xv^{2})^{2m+1}}$$ such
that $p_{m}^{\prime }(v_{0};x)$ is a power series, where $v_{0}=C(x)=\frac{1-%
\sqrt{1-4x}}{2x}$.
\end{theorem}

\begin{proof}
We prove this result by induction on $m$. For $m=0$, Theorem \ref{ppm2} for $%
m=0$ gives that $$\left( 1+\frac{xv^{2}}{1-v}\right)
p_{0}(v;x)=x^{3}+\frac{x}{1-v}p_{0}(1;x).
$$
Again, this type of equation can be solved with the kernel method.
Substituting $v=C(x)=\frac{1-\sqrt{1-4x}}{2x}$ in the above
equation, we get that $p_{0}(1;x)=x^{3}C^{2}(x)$. That is the number
of $0$-parallel permutations on $n$ letters is given by
$\frac{1}{n-1}{{2n-4}\choose{n-2}}$, for all $n\geq 3$. Moreover,
the ordinary generating function $p_{0}(v;x)$ is given by
$$
p_{0}(v;x)=\frac{(1-v)x^{3}+x^{4}C^{2}(x)}{1-v+xv^{2}},
$$
hence the theorem holds for $m=0$. Let $p_{m}(v;x)=\frac{p_{m}^{\prime }(v;x)%
}{(1-v+xv^{2})^{2m+1}}$, where $p_{m}^{\prime }(C(x);x)$ is a power
series. Then Theorem \ref{ppm2} gives that
$$
\left( 1+\frac{xv^{2}}{1-v}\right) p_{m+1}(v;x)=v^{2}x\frac{d}{dv}p_{m}(v;x)+%
\frac{x}{1-v}p_{m+1}(1;x),
$$
which is equivalent to
$$
\begin{array}{l}
\left( 1-v+xv^{2}\right) p_{m+1}(v;x)=v^{2}(1-v)x\frac{\frac{d}{dv}p_{m}^{\prime }(v;x)}{(1-v+xv^{2})^{2m+1}}\\
\qquad\qquad\qquad\qquad\qquad\quad-(2m+1)v^{2}(1-v)(1-2vx)x\frac{p_{m}^{\prime }(v;x)}{%
(1-v+xv^{2})^{2m+2}}+xp_{m+1}(1;x).\end{array}$$ Multiplying by
$(1-v+xv^{2})^{2m+2}$, we get that
$$
\begin{array}{l}
\left( 1-v+xv^{2}\right) ^{2m+3}p_{m+1}(v;x) \\
\qquad =v^{2}(1-v)(1-v+xv^{2})x\frac{d}{dv}p_{m}^{\prime
}(v;x)-(2m+1)v^{2}(1-v)(1-2vx)xp_{m}^{\prime }(v;x) \\
\qquad +x(1-v+xv^{2})^{2m+2}p_{m+1}(1;x).%
\end{array}%
$$
Now, differentiating the above recurrence relation $2m+2$ times respect to $%
v $, we can write
$$
\begin{array}{l}
\frac{d^{2m+2}}{dv^{2m+2}}\left[ \left( 1-v+xv^{2}\right)
^{2m+3}p_{m+1}(v;x)\right] \\
=\frac{d^{2m+2}}{dv^{2m+2}}\left[v^{2}(1-v)(1-v+xv^{2})x\frac{d}{dv}
p_{m}^{\prime }(v;x)-(2m+1)v^{2}(1-v)(1-2vx)xp_{m}^{\prime }(v;x)\right] \\
 +x\frac{d^{2m+2}}{dv^{2m+2}}\left[ (1-v+xv^{2})^{2m+2}\right]
p_{m+1}(1;x).
\end{array}%
$$
Substituting $v=C(x)$,
$$
{\small\begin{array}{l}
-(2m+2)!x(1-2xC(x))^{2m+2}p_{m+1}(1;x) \\
=\frac{d^{2m+2}}{dv^{2m+2}}\left[
v^{2}(1-v)(1-v+xv^{2})x\frac{d}{dv} p_{m}^{\prime
}(v;x)-(2m+1)v^{2}(1-v)(1-2vx)xp_{m}^{\prime }(v;x)\right]
\biggr|_{v=C(x).}
\end{array}}
$$
Since $p_{m}^{\prime }(v;x)$ is a generating function defined at
$v=C(x)$ then any derivative of $p_{m}^{\prime }(v;x)$ respect to
$v$ is a generating
function defined at $v=C(x)$, which gives an explicit formula for $%
p_{m+1}(1;x)$. If we substitute the formulas of $p_{m+1}(1;x)$ and $%
p_{m}(v;x)=\frac{p_{m}^{\prime }(v;x)}{(1-v+xv^{2})^{2m+1}}$ in the
functional equation
$$
\left( 1+\frac{xv^{2}}{1-v}\right) p_{m+1}(v;x)=v^{2}x\frac{d}{dv}p_{m}(v;x)+%
\frac{x}{1-v}p_{m+1}(1;x)
$$
we see that the generating function $p_{m+1}(v;x)$ can be written as $\frac{%
p_{m+1}^{\prime }(v;x)}{(1-v+xv^{2})^{2m+3}}$ such that
$p_{m+1}^{\prime
}(C(x);x)$ is a power series. Hence, the theorem is proved by induction on $%
m $.
\end{proof}

Theorem \ref{ppm3} provides an algorithm for finding $p_{m}(v;x)$ for any
given $m\geq 0$. From the proof of Theorem~\ref{ppm3}, with the help of any
scientific computing software, we can stated the following result.

\begin{theorem}
\label{ppm4} For $m=0,1,2,3,4,5$ the ordinary generating function $%
p_{m}(1;x) $ is given by
$$
\begin{array}{ll}
p_{0}(1;x)=\frac{x(1-2x)}{2}-\frac{x}{2}\sqrt{1-4x} &  \\[5pt]
p_{1}(1;x)=\frac{x(2x^{2}-4x+1)}{2(1-4x)}+\frac{x(2x-1)}{2\sqrt{1-4x}} &  \\%
[5pt]
p_{2}(1;x)=\frac{x(1-2x)(4x^{2}-6x+1)}{2(1-4x)^{2}}-\frac{%
x(6x^{4}-28x^{3}+30x^{2}-10x+1)}{2\sqrt{1-4x}^{5}} &  \\[5pt]
p_{3}(1;x)=\frac{x(24x^{6}-152x^{5}+300x^{4}-256x^{3}+96x^{2}-16x+1)}{%
2(1-4x)^{4}}+\frac{x(2x-1)(18x^{4}-48x^{3}+46x^{2}-12x+1)}{2\sqrt{1-4x}^{7}}
&  \\[5pt]
p_{4}(1;x)=\frac{x(1-2x)(96x^{6}-328x^{5}+496x^{4}-392x^{3}+124x^{2}-18x+1)}{%
2(1-4x)^{5}} &  \\
\qquad \qquad \qquad \qquad \qquad -\frac{%
x(126x^{8}-888x^{7}+2268x^{6}-3068x^{5}+2310x^{4}-924x^{3}+198x^{2}-22x+1)}{2%
\sqrt{1-4x}^{11}} &  \\[5pt]
p_{5}(1;x)=\frac{%
x(864x^{10}-6048x^{9}+17264x^{8}-28736x^{7}+30984x^{6}-21504x^{5}+8960x^{4}-2240x^{3}+336x^{2}-28x+1)%
}{2(1-4x)^{7}} &  \\
\qquad \qquad \qquad \qquad \quad +\frac{%
x(2x-1)(630x^{8}-2352x^{7}+4404x^{6}-4960x^{5}+3526x^{4}-1240x^{3}+238x^{2}-24x+1)%
}{2\sqrt{1-4x}^{13}} &
\end{array}%
$$
\end{theorem}

We remark that it is not hard to prove by induction, as the proof of Theorem %
\ref{ppm3}, that our generating function $p_{m}(1;x)$ is a rational function
in the variables $x$ and $\sqrt{1-4x}$.

\section{Enumeration of three faces polygons}

A permutation $\pi $ is said to be \emph{triangular} if it begins at letter $%
1$ and its polygon $P_{\pi }$ has at most $3$ faces. For example, there
exist $1,1,6,20$ triangular permutations of length $1,2,3,4$, respectively.
We denote the set of all triangular permutations of length $n$ by $_{n}$.
Given $a_{1},a_{2},\ldots ,a_{d}\in \mathbb{N}$, we define 
\[
t_{n;a_{1},a_{2},\ldots ,a_{d}}=\#\{\pi _{1}\pi _{2}\ldots
\pi_{n}\in_{n}\mid \pi _{1}\pi _{2}\ldots \pi _{d}=a_{1}a_{2}\ldots a_{d}\}, 
\]
The cardinality of the set $_{n}$ is denoted by $t_{n}$.

\begin{theorem}
\label{gg1} The number of triangular permutations of length $n+2$ is
${{ 2n}\choose{n}}$. Moreover, the ordinary generating function
$t(v;x)=\sum_{n\geq 2}\sum_{a=2}^{n}t_{n;1,a}v^{a-2}x^{n}$ is given
by
$$
\frac{x^{2}(1-v)(1-xv)^{2}}{(1-2xv)(1-v+xv^{2})}+\frac{x^{3}}{1-v+xv^{2}}%
\cdot \frac{1}{\sqrt{1-4x}}.
$$
\end{theorem}

\begin{proof}
From the definitions, we have that $
t_{n}=t_{n;1}=t_{n;1,2}+t_{n;1,3}+\cdots +t_{n;1,n}$. For all
$a=3,4,\ldots ,n-1$,
\begin{equation}
t_{n;1,a}=t_{n;1,a,2}+\sum_{j=a+1}^{n}t_{n;1,a,j}=t_{n-1;1,a-1}+%
\sum_{j=a+1}^{n}t_{n-1;1,j-1}=\sum_{j=a-1}^{n-1}t_{n-1;1,j},
\label{eqgg22}
\end{equation}%
with the initial conditions $t_{n;1,2}=t_{n-1;1}$ and
$t_{n;1,n}=2^{n-3}$. To see that
\begin{equation}
t_{n;1,n}=2^{n-3}  \label{eqgg1}
\end{equation}%
we consider the following equation
$t_{n;1,n}=t_{n;1,n,2}+t_{n;1,n,n-1}=t_{n-1;1,n-1}+t_{n-1;1,n-1}=2t_{n-1;1,n-1}$
for all $n\geq 4$, and $t_{3;1,3}=1$ which implies that
$t_{n;1,n}=2^{n-3}$ as claimed.

To solve the recurrence relation of the sequence $t_{n;1,a}$, we
need to define $t_{n}(v)=\sum_{a=2}^{n}t_{n;1,a}v^{a-2}$. Thus,
multiplying the above recurrence relation by $v^{a-2}$ and summing
over $a=3,4,\ldots ,n-1$ we obtain that
$$
t_{n}(v)=t_{n-1}(1)+2^{n-3}v^{n-2}+\sum_{a=3}^{n-1}%
\sum_{j=a-1}^{n-1}t_{n-1;1,j}v^{a-2},
$$%
which is equivalent to
$$
t_{n}(v)=t_{n-1}(1)+2^{n-4}v^{n-2}+\frac{v}{1-v}(t_{n-1}(1)-vt_{n-1}(v)),
$$%
for $n\geq 4$. Let $t(v;x)=\sum_{n\geq 2}t_{n}(v)x^{n}$. If
multiplying the above recurrence relation with $x^{n}$ and summing
over all possibly $n\geq 4 $ by using the initial conditions
$t_{2}(v)=1$ and $t_{3}(v)=1+v$, we then obtain the following
functional equation
$$
t(v;x)=\frac{xv}{1-v}(t(1;x)-vt(v;x))+xt(1;x)-\frac{(1-xv)^{2}x^{2}}{1-2xv},
$$%
which is equivalent to
$$
\left( 1+\frac{xv^{2}}{1-v}\right) t(v;x)=\frac{x}{1-v}t(1;x)-\frac{%
(1-xv)^{2}x^{2}}{1-2xv}.
$$%
This type of equation can be solved using the kernel method. Substitute $v=%
\frac{1-\sqrt{1-4x}}{2x}$ in the above functional equation to get $t(1;x)=%
\frac{x^{2}}{\sqrt{1-4x}}$, that is, the number of triangular
permutations of length $n$ is exactly ${{2n-4}\choose{n-2}}$, as
required. Moreover, substituting the expression of $t(1;x)$ in the
functional equation, we get an explicit formula for $t(v;x)$, as
claimed.
\end{proof}

As a corollary of Theorem \ref{pp1} and Theorem\ref{gg1} we get the
following.

\begin{corollary}
\label{co33} The number of polygons on $n$ vertices with exactly
three faces is $\frac{4(n-2)}{n-1}{{2n-4}\choose{n-2}}$, for all
$n\geq2$.
\end{corollary}

\begin{proof}
Theorem \ref{pp1} and Theorem \ref{gg1} give that the number of
permutations $\pi $ of length $n$ that begin at letter $1$ and its
polygon $P_{\pi }$ has exactly three faces is
$\frac{n-2}{n-1}{{2n-4}\choose{n-2}}$. If $P_{\pi }$ has exactly
three faces then also $P_{\pi ^{\prime }}$ and $P_{\pi ^{\prime
\prime }}$ have exactly three faces, where $\pi ^{\prime }$ is the
complement of $\pi $ and $\pi ^{\prime \prime }$ is the reversal of
$\pi $. (Recall that the \emph{reversal} of a permutation $\pi
_{1}\pi _{2}\ldots \pi _{n}$ is $\pi _{n}\ldots \pi _{2}\pi _{1}$;
the \emph{complement} of is the permutation $(n+1-\pi _{1})(n+1-\pi
_{2})\ldots (n+1-\pi _{n})$). From this fact, we obtain that the
number of polygons on $n$ vertices with exactly three faces is four
times the number of permutations $\pi $ of length $n$ that begin at
letter $1$ and whose polygon $P_{\pi }$ has exactly three faces.
This number is $\frac{4(n-2)}{n-1}{{2n-4}\choose{n-2}}$, as
required.
\end{proof}

Let us now generalize the above enumeration of polygons with exactly three
faces. A polygon is said to be \emph{$m$-isolated polygon} if its a grid
polygon contains exactly $m$ isolated points. A permutation $\pi $ is said
to be \emph{$m$-triangular} if it begins at letter $1$ and by deleting $m$
vertices from $L_{\pi }$ gives a polygon with exactly three faces. Clearly a 
$0$-triangular permutation is a triangular permutation as defined above. In
Table \ref{tab2} we give the number of $m$-triangular permutations of length 
$n$ for $m\leq 5$ and $n\leq 13$. 
\begin{table}[h]
\begin{center}
{\small 
\begin{tabular}{c|lllllllllllll}
$m/n$ & $0$ & $1$ & $2$ & $3$ & $4$ & $5$ & $6$ & $7$ & $8$ & $9$ & $10$ & $%
11$ & $12$ \\ \hline
$0$ & $1$ & $1$ & $1$ & $2$ & $6$ & $20$ & $70$ & $252$ & $924$ & $3432$ & $%
12870$ & $48620$ & $184756$ \\ 
$1$ & $0$ & $0$ & $0$ & $0$ & $0$ & $4$ & $42$ & $300$ & $1812$ & $9960$ & $%
51546$ & $255868$ & $1231932$ \\ 
$2$ & $0$ & $0$ & $0$ & $0$ & $0$ & $0$ & $8$ & $144$ & $1572$ & $13440$ & $%
99042$ & $660068$ & $4090848$ \\ 
$3$ & $0$ & $0$ & $0$ & $0$ & $0$ & $0$ & $0$ & $24$ & $636$ & $9576$ & $%
107718$ & $1007884$ & $8295360$ \\ 
$4$ & $0$ & $0$ & $0$ & $0$ & $0$ & $0$ & $0$ & $0$ & $96$ & $3432$ & $66936$
& $945152$ & $10827036$ \\ 
$5$ & $0$ & $0$ & $0$ & $0$ & $0$ & $0$ & $0$ & $0$ & $0$ & $480$ & $21888$
& $529912$ & $9084708$%
\end{tabular}%
}
\end{center}
\caption{Number of $m$-triangular permutations of length $n$.}
\label{tab2}
\end{table}

We denote the set of all $m$-triangular permutations of length $n$ by $%
\mathcal{T}_{n}^{m}$. Given $a_{1},a_{2},\ldots ,a_{d}\in \mathbb{N}$, we
define 
\[
r_{n,m;a_{1},a_{2},\ldots ,a_{d}}=\#\{\pi _{1}\pi _{2}\ldots \pi _{n}\in%
\mathcal{T}_{n}^{m}\mid \pi _{1}\pi _{2}\ldots \pi _{d}=a_{1}a_{2}\ldots
a_{d}\}, 
\]%
The cardinality of the set $\mathcal{T}_{n}^{m}$ is defined by $r_{n,m}$.

\begin{theorem}
\label{ggm1} Let $r(v;q,x)=\sum_{n\geq
4}\sum_{a=2}^{n}\sum_{m=0}^{n}r_{n,m;1,a}q^{m}v^{a-2}x^{n}$ be the
ordinary generating function for the sequence $r_{n,m;1,a}$. Then
\begin{eqnarray*}
&&\left( 1+\frac{xv^{2}}{1-v}\right) r(v;q,x)-v^{2}xq\frac{d}{dv}r(v;q,x) \\
&&\qquad\qquad=2\frac{1-v^{3}}{1-v}x^{4}+\frac{x}{1-v}r(1;q,x)-(1-q)x^{2}\sum_{n\geq
3}(vx)^{n}\prod_{j=0}^{n-3}(2+jq).
\end{eqnarray*}
\end{theorem}

\begin{proof}
Define $r_{n;a_{1},a_{2},\ldots ,a_{d}}=\sum_{m\geq
0}r_{n,m;a_{1},\ldots ,a_{d}}q^{m}$ for any $n$ and $d$. From the
definitions we have that
$$
r_{n}=r_{n;1}=r_{n;1,2}+r_{n;1,3}+\cdots +r_{n;1,n}.
$$%
For all $a=3,4,\ldots ,n-1$,
$$
\begin{array}{ll}
r_{n;1,a} & =r_{n;1,a,2}+\sum\limits_{j=a+1}^{n}r_{n;1,a,j}+%
\sum\limits_{j=3}^{a-1}r_{n;1,a,j} \\
& =r_{n-1;1,a-1}+\sum\limits_{j=a+1}^{n}r_{n-1;1,j-1}+q(a-3)r_{n-1;1,a-1} =q(a-3)r_{n-1;1,a-1}+\sum\limits_{j=a-1}^{n-1}r_{n-1;1,j},%
\end{array}%
$$%
with the initial conditions $r_{n;1,2}=r_{n-1;1}$ and $%
r_{n;1,n}=(2+q(n-4))r_{n-1;1,n-1}$. To see the last quantity $%
r_{n;1,n}=(2+q(n-4))r_{n-1;1,n-1}$ we consider the following equation $%
r_{n;1,n}=r_{n;1,n,2}+r_{n;1,n,n-1}+\sum_{j=3}^{n-2}r_{n;1,n,j}$
which leads to
$r_{n;1,n}=r_{n-1;1,n-1}+r_{n-1;1,n-1}+q(n-4)r_{n-1;1,n-1}$, or
equivalently $r_{n;1,n}=(2+q(n-4))r_{n-1;1,n-1}$, for all $n\geq 4$, and $%
r_{3;1,3}=1$.

To solve the recurrence relation of the sequence $r_{n;1,a}$, we
need to define the generating function
$r_{n}(v;q)=\sum_{a=2}^{n}r_{n;1,a}v^{a-2}$. Thus, multiplying the
above recurrence relation by $v^{a-2}$ and summing over
$a=3,4,\ldots ,n-1$, for all $n\geq 5$, we obtain
$$\begin{array}{l}
r_{n}(v;q)=r_{n-1}(1)+\frac{v}{1-v}(r_{n-1}(1;q)-vr_{n-1}(v;q))+v^{2}q\frac{d%
}{dv}(r_{n-1}(v;q))+(1-q)v^{n-2}\prod\limits_{j=0}^{n-5}(2+jq),
\end{array}$$
with the initial condition $r_{4}(v;q)=2+2v+2v^{2}$. Let $%
r(v;q,x)=\sum_{n\geq 4}r_{n}(v;q)x^{n}$. If multiplying the above
recurrence relation with $x^{n}$ and summing over all $n\geq 5$, we
then obtain the requested functional equation.
\end{proof}

It is well known that the unsigned Stirling numbers $s_{n,k}$ of the first
kind (these count the number of ways to permute a list of $n$ items into $k$
cycles \cite[Sequence A008275]{JS}), satisfy $\prod_{j=0}^{n}(p+j)=%
\sum_{j=1}^{n+1}s_{n+1,j}p^{j}$. This identity leads to that 
\[
\sum_{n\geq 3}(vx)^{n}\prod_{j=0}^{n-3}(2+jq)=\sum_{n\geq
3}(vx)^{n}\sum_{j=1}^{n-2}s_{n-2,j}2^{j}q^{n-2-j}. 
\]%
Theorem \ref{ggm1} gives a recurrence relation for the ordinary generating
function $r_{m}(v;x)$ for the number of $m$-triangular permutations of
length $n$, as follows.

\begin{theorem}
\label{ggm2} Let $r_{m}(v;x)$ be the coefficient of $q^{m}$ in the
ordinary generating function $r(v;q,x)$, that is,
$r_{m}(v;x)=\sum_{n\geq
4}\sum_{a=2}^{n}r_{n,m;1,a}q^{m}v^{a-2}x^{n}$ (define
$r_{-1}(v;x)=0$). Then
\begin{eqnarray*}
\left( 1+\frac{xv^{2}}{1-v}\right) r_{m}(v;x) &=&v^{2}x\frac{d}{dv}%
r_{m-1}(v;x)+2\delta _{m,0}\frac{1-v^{3}}{1-v}x^{4} \\
&&+\frac{x}{1-v}r_{m}(1;x)+x^{2}(d_{m}(vx)-d_{m-1}(vx)),
\end{eqnarray*}%
where $\delta _{0,0}=1$, $\delta _{m,0}=0$ with $m\neq0$, and
$d_{m}(y)=\sum_{n\geq 3}2^{n-2-m}y^{n}s_{n-2,n-2-m}$ where $s_{n,k}$
are the unsigned Stirling numbers of the first kind.
\end{theorem}

Theorem \ref{ggm2} can be used to obtain an explicit formula for $r_{m}(v;x)$
for given $m$. First, we need the following lemma.

\begin{lemma}
\label{gglm1} For all $n\geq 0$,
$$
\begin{array}{ll}
s_{n,n}=1, & s_{n,n-3}=\frac{n^{2}(n-1)^{2}(n-2)(n-3)}{48} \\
s_{n,n-1}=\frac{n(n-1)}{2}, & s_{n,n-4}=\frac{%
n(n-1)(n-2)(n-3)(n-4)(15n^{3}-30n^{2}+5n+2)}{5760} \\
s_{n,n-2}=\frac{n(n-1)(n-2)(3n-1)}{24}, & s_{n,n-5}=\frac{%
n^{2}(n-1)^{2}(n-2)(n-3)(n-4)(n-5)(3n^{2}-7n-2)}{11520}.%
\end{array}%
$$%
Moreover, the ordinary generating function $d_{k}(x)=\sum_{n\geq
3}2^{n-2-k}x^{n}s_{n-2,n-2-k}$ is a rational function with only one pole at $%
x=\frac{1}{2}$.
\end{lemma}

\begin{proof}
It is well known that the unsigned Stirling numbers $s_{n,k}$ of the
first kind satisfy the recurrence relation
$s_{n,k}=s_{n-1,k-1}+(n-1)s_{n-1,k}$,
for all $k=1,2,\ldots ,n-1$ and with the initial conditions $%
s_{n,n}=s_{n,0}=1$. We the help of any scientific computing software
one can obtain the requested result. Using the recurrence relation
and induction, we
get that the ordinary generating function $d_{k}(x)$ has only one pole at $x=%
\frac{1}{2}$, as required.
\end{proof}

Theorem \ref{ggm2} provides an algorithm for finding $r_{m}(v;x)$ for any
given $m\geq 0$, since we are dealing with a functional equation with one
variable, and this type of functional equations can be solved using the
kernel method. We remark that we cannot just substitute $v=C(x)=\frac{1-%
\sqrt{1-4x}}{2x}$ in the functional equation of the statement of Theorem \ref%
{ggm2}, since it may well be that the generating function $\frac{d}{dv}%
r_{m-1}(v;x)$ is not defined at $v=C(x)$. This kind of recurrence relation
for given $m$ can be solved with the following result.

\begin{theorem}
\label{ggm3} For any $m\geq 0$, the ordinary generating function
$r_{m}(v;x)$ can be written as $$\frac{r_{m}^{\prime
}(v;x)}{(1-v+xv^{2})^{2m+1}}$$ such
that $r_{m}^{\prime }(v_{0};x)$ is a power series, where $v_{0}=C(x)=\frac{1-%
\sqrt{1-4x}}{2x}$.
\end{theorem}

\begin{proof}
We prove the theorem by induction on $m$. For $m=0$, Theorem
\ref{ggm2} together with Lemma \ref{gglm1} for $m=0$ give that
$$
\left( 1+\frac{xv^{2}}{1-v}\right) r_{0}(v;x)=2\frac{1-v^{3}}{1-v}x^{4}+%
\frac{2v^{3}x^{5}}{2vx-1}.
$$%
This type of equation can be solved using the kernel method. Substitute $%
v=C(x)=\frac{1-\sqrt{1-4x}}{2x}$ in the above functional equation to
get that
$$
r_{0}(1;x)=\frac{x^{2}}{\sqrt{1-4x}}-x^{2}(2x+1).
$$%
Then, the number of $0$-triangular permutations on $n$ letters is
given by ${{2n-4}\choose{n-2}}$, for all $n\geq 2$. Moreover, the
ordinary generating function $r_{0}(v;x)$ is given by
$$
r_{0}(v;x)=\frac{x^{3}(1-2vx)+x^{3}(2v^{4}x^{2}+2v^{3}x^{2}-2xv^{3}+2xv-1)%
\sqrt{1-4x}}{(1-2vx)(1-v+xv^{2})\sqrt{1-4x}},
$$%
hence the theorem holds for $m=0$. Let $r_{m}(v;x)=\frac{r_{m}^{\prime }(v;x)%
}{(1-v+xv^{2})^{2m+1}}$, where $r_{m}^{\prime }(C(x);x)$ is a power
series. Then Theorem \ref{ggm2} gives that
$$
\left( 1+\frac{xv^{2}}{1-v}\right) r_{m+1}(v;x)=v^{2}x\frac{d}{dv}r_{m}(v;x)+%
\frac{x}{1-v}r_{m+1}(1;x)+x^{2}(d_{m+1}(vx)-d_{m}(vx)),
$$%
which is equivalent to
$$
\begin{array}{l}
\left( 1-v+xv^{2}\right) r_{m+1}(v;x) \\
\qquad =v^{2}(1-v)x\frac{\frac{d}{dv}r_{m}^{\prime }(v;x)}{%
(1-v+xv^{2})^{2m+1}}-(2m+1)v^{2}(1-v)(1-2vx)x\frac{r_{m}^{\prime }(v;x)}{%
(1-v+xv^{2})^{2m+2}} \\
\qquad +xr_{m+1}(1;x)+x^{2}(1-v)(d_{m+1}(vx)-d_{m}(vx)).%
\end{array}%
$$%
Multiplying by $(1-v+xv^{2})^{2m+2}$, we obtain
$$
\begin{array}{l}
\left( 1-v+xv^{2}\right) ^{2m+3}r_{m+1}(v;x) \\
=v^{2}(1-v)(1-v+xv^{2})x\frac{d}{dv}r_{m}^{\prime
}(v;x)-(2m+1)v^{2}(1-v)(1-2vx)xr_{m}^{\prime }(v;x) \\
+x(1-v+xv^{2})^{2m+2}r_{m+1}(1;x)+x^{2}(1-v)(1-v+xv^{2})^{2m+2}(d_{m+1}(vx)-d_{m}(vx)).%
\end{array}%
$$%
Now, differentiating the above recurrence relation $2m+2$ times respect to $%
v $ we get that
$$
\begin{array}{l}
\frac{d^{2m+2}}{dv^{2m+2}}\left[ \left( 1-v+xv^{2}\right) ^{2m+3}r_{m+1}(v;x)%
\right] \\
=\frac{d^{2m+2}}{dv^{2m+2}}\left[ v^{2}(1-v)(1-v+xv^{2})x\frac{d}{dv}%
r_{m}^{\prime }(v;x)-(2m+1)v^{2}(1-v)(1-2vx)xr_{m}^{\prime }(v;x)\right] \\
+x\frac{d^{2m+2}}{dv^{2m+2}}\left[ (1-v+xv^{2})^{2m+2}\right]
r_{m+1}(1;x) \\
+\frac{d^{2m+2}}{dv^{2m+2}}\left[
x^{2}(1-v)(1-v+xv^{2})^{2m+2}(d_{m+1}(vx)-d_{m}(vx))\right] .%
\end{array}%
$$%
Substituting $v=C(x)$,
$$
{\small\begin{array}{l}
-(2m+2)!x(1-2xC(x))^{2m+2}r_{m+1}(1;x) \\
=\frac{d^{2m+2}}{dv^{2m+2}}\left[ v^{2}(1-v)(1-v+xv^{2})x\frac{d}{dv}%
r_{m}^{\prime }(v;x)-(2m+1)v^{2}(1-v)(1-2vx)xr_{m}^{\prime }(v;x)\right] %
\biggr|_{v=C(x)} \\
+\frac{d^{2m+2}}{dv^{2m+2}}\left[
x^{2}(1-v)(1-v+xv^{2})^{2m+2}(d_{m+1}(vx)-d_{m}(vx))\right] \biggr|_{v=C(x).}%
\end{array}}
$$%
Since $r_{m}^{\prime }(v;x)$ is a generating function defined at
$v=C(x)$ then any derivative of $r_{m}^{\prime }(v;x)$ respect to
$v$ is a generating function defined at $v=C(x)$. Lemma \ref{gglm1}
provides that the above recurrence relation gives to an explicit
formula for $r_{m+1}(1;x)$. If we
now substitute the formula of $r_{m+1}(1;x)$ and $r_{m}(v;x)=\frac{%
r_{m}^{\prime }(v;x)}{(1-v+xv^{2})^{2m+1}}$ in the functional
equation
$$
\left( 1+\frac{xv^{2}}{1-v}\right) r_{m+1}(v;x)=v^{2}x\frac{d}{dv}r_{m}(v;x)+%
\frac{x}{1-v}r_{m+1}(1;x)+x^{2}(d_{m+1}(vx)-d_{m}(vx)),
$$%
we obtain that the generating function $r_{m+1}(v;x)$ can be written as $%
\frac{r_{m+1}^{\prime }(v;x)}{(1-v+xv^{2})^{2m+3}}$, such that $%
r_{m+1}^{\prime }(C(x);x)$ is power series. Hence, the theorem is
proved by induction on $m$.
\end{proof}

In conjunction with the kernel method, Theorem \ref{ggm3} provides an
algorithm for finding $r_{m}(v;x)$ for any given $m\geq 0$. From the proof
of Theorem \ref{ggm3}, using any scientific computing software, we can state
the following.

\begin{theorem}
\label{ggm4} For $m=0,1,2,3,4,5$ the ordinary generating function $%
r_{m}(1;x) $ is given by
$$
\begin{array}{ll}
r_{0}(1;x)=-x^{2}(2x+1)+\frac{x^{2}}{\sqrt{1-4x}} &  \\[5pt]
r_{1}(1;x)=\frac{-x^{2}(2x^{2}+2x-1)}{(1-4x)^{2}}-\frac{1}{\sqrt{1-4x}^{3}}
&  \\[5pt]
r_{2}(1;x)=\frac{x^{2}(4x^{3}+12x^{2}-6x-1)}{(1-4x)^{3}}+\frac{%
x^{2}(6x^{4}+4x^{3}-26x^{2}+4x+1)}{\sqrt{1-4x}^{7}} &  \\[5pt]
r_{3}(1;x)=\frac{-x^{2}(24x^{6}-56x^{5}-188x^{4}+136x^{3}+68x^{2}-22x-1)}{%
(1-4x)^{5}}-\frac{x^{2}(24x^{5}+18x^{4}-120x^{3}-18x^{2}+24x+1)}{\sqrt{1-4x}%
^{9}} &  \\[5pt]
r_{4}(1;x)=\frac{%
x^{2}(144x^{7}-208x^{6}-720x^{5}+360x^{4}+1064x^{3}-204x^{2}-66x-1)}{%
(1-4x)^{6}} &  \\
\qquad \qquad \qquad \qquad \qquad +\frac{%
x^{2}(126x^{8}-504x^{7}-132x^{6}+2064x^{5}+1086x^{4}-1608x^{3}+70x^{2}+64x+1)%
}{(1-4x)^{13}} &  \\[5pt]
r_{5}(1;x)=\frac{%
-x^{2}(+864x^{10}-4128x^{9}+4688x^{8}+9536x^{7}-3304x^{6}-42736x^{5}+12584x^{4}+7304x^{3}-1444x^{2}-154x-1)%
}{(1-4x)^{8}} &  \\
\qquad \qquad \qquad \qquad\qquad  -\frac{%
x^{2}(1008x^{9}-2538x^{8}-336x^{7}+8004x^{6}+13656x^{5}-15378x^{4}-3472x^{3}+1758x^{2}+156x+1)%
}{\sqrt{1-4x}^{15}} &
\end{array}%
$$
\end{theorem}

We remark that it is not hard to prove by induction, as the proof of Theorem %
\ref{ggm3}, that our generating function $p_{m}(1;x)$ is a rational function
in the variables $x$ and $\sqrt{1-4x}$.

\section{Enumeration of four faces polygons}

A permutation $\pi $ is said to be \emph{square} if the subsequence of the
sources of $L_{\pi }$ lies on at most two faces of $P_{\pi }$. For example,
there exists $1,2,6,24,104,464,2088$ square permutations of length $%
1,2,3,4,5,6,7$, respectively. We denote the set of all square permutations
of length $n$ by $\mathcal{Q}_{n}$. Given $a_{1},a_{2},\ldots ,a_{d}\in 
\mathbb{N}$, we define 
\[
q_{n;a_{1},a_{2},\ldots ,a_{d}}=\#\{\pi _{1}\pi _{2}\ldots \pi _{n}\in 
\mathcal{Q}_{n}\mid \pi _{1}\pi _{2}\ldots \pi _{d}=a_{1}a_{2}\ldots
a_{d}\}, 
\]
The cardinality of the set set $\mathcal{Q}_{n}$ by $q_{n}$. Clearly, a
triangular permutation is a square permutation. We derive an explicit
formula for the number of square permutations of length $n$ as follows.

\begin{theorem}
\label{hhs4} The ordinary generating function for the number of
square permutations of length $n$ is given by
$$
1+x+\frac{2(1-3x)x^{2}}{(1-4x)^{2}}-\frac{4x^{3}}{(1-4x)^{3/2}}.
$$%
Moreover, the number of square permutation of length $n$ is
$$
2(n+2)4^{n-3}-4(2n-5){{2n-6}\choose{n-3}},
$$%
for all $n\geq 3$.
\end{theorem}

\begin{proof}
From the symmetry arising in the construction of square permutations
we have that for all $n\geq a>b\geq 1$,
\begin{equation}
q_{n;a,b}=q_{n;n+1-a,n+1-b}\mbox{ and }q_{n;a,b}=q_{n;b,a}.
\label{eqff1}
\end{equation}%
Define
$Q_{n}(u,v)=\sum_{a=1}^{n}\sum_{b=1}^{n}q_{n;a,b}v^{a-1}u^{b-1}$,
for all $n\geq 2$, and $Q(u,v;x)=\sum_{n\geq 0}Q_{n}(v,u)x^{n}$ to
be the
ordinary generating function for the sequence $Q_{n}(u,v)$. Thus, (\ref%
{eqff1}) gives
\begin{equation}
Q(v,u;x)=Q^{\prime }(v,u;x)+Q^{\prime }(u,v;x),  \label{eqff2}
\end{equation}%
where
$$
Q^{\prime }(u,v;x)=\sum_{n\geq 2}Q_{n}^{\prime
}(v,u)x^{n}=\sum_{n\geq
2}x^{n}\sum_{a=2}^{n}\sum_{b=1}^{a-1}q_{n;a,b}v^{a-1}u^{b-1}.
$$%
To find an explicit formula for $Q^{\prime }(1,1;x)$, which leads us
to explicit formula for $Q(1,1,;x)$, the ordinary generating
function for the number of square permutations of length $n$, we
need to divide the generating function $Q^{\prime }(u,v;x)$ into
three parts. For all $n\geq a>b\geq 1$, define
$$
\begin{array}{ll}
A(v;x)=\sum\limits_{n\geq 2}A_{n}(v)x^{n} & =\sum\limits_{n\geq
2}x^{n}\sum\limits_{a=2}^{n}q_{n;a,1}v^{a-1}, \\
B(v;x)=\sum\limits_{n\geq 2}B_{n}(v)x^{n} & =\sum\limits_{n\geq
3}x^{n}\sum\limits_{b=2}^{n-1}q_{n;n,b}v^{b-1}, \\
C(v,u;x)=\sum\limits_{n\geq 2}C_{n}(u,v)x^{n} & =\sum\limits_{n\geq
4}x^{n}\sum\limits_{a=3}^{n-1}\sum\limits_{b=2}^{i-1}q_{n;a,b}v^{a-1}u^{b-1}.%
\end{array}%
$$%
Clearly, for all $n\geq 2$, $Q_{n}^{\prime
}(v,u)=C_{n}(v,u)+v^{n-1}B_{n}(u)+A_{n}(v)$ and then
\begin{equation}
Q^{\prime }(u,v;x)=C(v,u;x)+\frac{1}{v}B(u;xv)+A(v;x).
\label{eqff3}
\end{equation}

\textbf{Expression for $A(v;x)$}: First, we find an explicit formula
for the
ordinary generating function $A(v;x)$. From the definitions and (\ref{eqff1}%
), we have that
$$
\begin{array}{l}
q_{n;2,1}=q_{n-1;1}=\sum\limits_{b=2}^{n-1}q_{n-1;1,b}=%
\sum\limits_{b=2}^{n-1}q_{n-1;b,1}=A_{n-1}(1), \\
q_{n;a,1}=q_{n;a,1,2}+\sum\limits_{b=a+1}^{n}q_{n;a,1,b}=q_{n-1;a-1,1}+%
\sum\limits_{b=a+1}^{n}q_{n-1;1,b-1}=\sum\limits_{b=a-1}^{n-1}q_{n-1;b,1}, \\
q_{n;n,1}=q_{n;n,1,2}+q_{n;n,1,n-1}=2q_{n-1;n-1,1}.%
\end{array}%
$$%
Using $q_{3;3,1}=1$ and the recurrence relation for the sequence $q_{n;n,1}$%
, we obtain that, for all $n\geq 3$,
\begin{equation}
q_{n;n,1}=2^{n-3}.  \label{eqff4}
\end{equation}%
Multiplying by $v^{a-1}$ and summing over all $a=3,4,\ldots ,n-1$,
we obtain that
$$
\begin{array}{ll}
A_{n}(v) & =vA_{n-1}(1)+\sum\limits_{a=3}^{n-1}v^{a-1}%
\sum\limits_{j=a-1}^{n-1}q_{n-1;j,1}+q_{n;n,1}v^{n-1}, \\
& =vA_{n-1}(1)+\sum\limits_{a=2}^{n-2}q_{n-1;a,1}\frac{v^{2}-v^{i+1}}{1-v}+\frac{%
v^{2}-v^{n-1}}{1-v}q_{n-1;n-1,1}+q_{n;n,1}v^{n-1}.%
\end{array}%
$$%
Then (\ref{eqff4}) leads us to
$$\begin{array}{l}
A_{n}(v)=vA_{n-1}(1)+\frac{v^{2}}{1-v}(A_{n-1}(1)-2^{n-4}-A_{n-1}(v)+2^{n-4}v^{n-2})+2^{n-4}\frac{v^{2}-v^{n-1}}{1-v}+2^{n-3}v^{n-1},\end{array}
$$%
which is equivalent to
$$
A_{n}(v)=vA_{n-1}(v)+\frac{v^{2}}{1-v}(A_{n-1}(1)-A_{n-1}(v))+2^{n-4}v^{n-4},
$$%
for all $n\geq 4$, with initial conditions $A_{2}(v)=v$ and $%
A_{3}(v)=v+v^{2} $. Writing the above recurrence relation in terms
of generating functions,
$$
A(v;x)-(v+v^{2})x^{3}-vx^{2}=vx(A(1;x)-x^{2})+\frac{xv^{2}}{1-v}%
(A(1;x)-x^{2}-A(v;x)+vx^{2})+\frac{v^{3}x^{4}}{1-2vx}.
$$%
Equivalently,
$$
\left( 1+\frac{v^{2}x}{1-v}\right) A(v;x)=vx^{2}+\frac{v^{3}x^{4}}{1-2vx}+%
\frac{vx}{1-v}A(1;x).
$$%
This type of equation can be solved systematically using the kernel
method. We substitute $v=\frac{1-\sqrt{1-4x}}{2x}$ in the above
functional equation to get $A(1;x)=\frac{x^{2}}{\sqrt{1-4x}}$ and
then
\begin{equation}
A(v;x)=\frac{1}{1-v+v^{2}x}\left( v(1-v)x^{2}+\frac{v^{3}(1-v)x^{4}}{1-2vx}+%
\frac{vx^{3}}{\sqrt{1-4x}}\right) .  \label{eqff5}
\end{equation}

\textbf{Expression for $B(v;x)$}: Using the symmetry on the set of
square permutations, see (\ref{eqff1}), we obtain that
$$
B(v;x)=\sum_{n\geq
3}x^{n}\sum_{j=2}^{n-1}q_{n;j,n}v^{j-1}=\sum_{n\geq
3}x^{n}\sum_{j=2}^{n-1}q_{n;n+1-j,1}v^{j-1}=\sum_{n\geq
3}x^{n}\sum_{j=2}^{n-1}q_{n;j,1}v^{n-j},
$$%
and from the definition of the generating function $A(v;x)$ together with (%
\ref{eqff4}),
$$
B(v;x)=\frac{1}{v}A(1/v;vx)-\frac{x^{2}(1-x)}{1-2x}.
$$%
It follows that
\begin{equation}
B(v;x)=\frac{vx^{3}(1-x)}{(1-2x)(1-v-x)}-\frac{x^{3}v^{2}}{(1-v-x)\sqrt{1-4vx%
}}.  \label{eqff6}
\end{equation}

\textbf{Expression for $C(v,u;x)$}: From the definitions and
(\ref{eqff1}), for all $n-1\geq a>b\geq 2$,
$$\begin{array}{ll}
q_{n;a,b}&=\sum\limits_{j=1}^{b-1}q_{n;a,b,j}+\sum\limits_{j=a+1}^{n}q_{n;a,b,j}\\
&=\sum\limits_{j=1}^{b-1}q_{n-1;a-1,j}+\sum\limits_{j=a+1}^{n}q_{n-1;j-1,b}=\sum\limits_{j=1}^{b-1}q_{n-1;a-1,j}+\sum\limits_{j=a}^{n-1}q_{n-1;j,b}.\end{array}
$$%
Thus, for all $n\geq 5$,
$$
\begin{array}{ll}
C_{n}(v,u) & =\sum\limits_{a=3}^{n-1}\sum\limits_{b=2}^{a-1}\left(
\sum\limits_{j=1}^{b-1}q_{n-1;a-1,j}+\sum\limits_{j=a}^{n-1}q_{n-1;j,b}\right)
v^{a-1}u^{b-1} \\
& =\sum\limits_{a=2}^{n-2}\sum\limits_{b=1}^{a-1}q_{n-1;a,b}\frac{u^{j}-u^{i}}{1-u}%
v^{i}+\sum\limits_{a=3}^{n-1}\sum\limits_{b=2}^{a-1}q_{n-1;a,b}\frac{v^{j}-v^{i}}{1-v}%
u^{j-1}.%
\end{array}%
$$%
Therefore, by the definition of the sequences $A_{n}(v)$, $B_{n}(v)$ and $%
C_{n}(n,u)$ together with (\ref{eqff1}), for all $n\geq 5$,
$$
\begin{array}{ll}
C_{n}(v,u) & =\frac{vu}{1-u}(C_{n-1}(v,u)-C_{n-1}(vu,1))+\frac{v}{1-v}%
(C_{n-1}(1,vu)-C_{n-1}(v,u)) \\
&\qquad+\frac{v}{1-v}(B_{n-1}(vu)-v^{n-2}B_{n-1}(u))+\frac{uv}{1-u}%
(A_{n-1}(v)-A_{n-1}(vu))\\
&\qquad-2^{n-4}\frac{uv^{n-1}(1-u^{n-2})}{1-u}.%
\end{array}%
$$%
By converting the above recurrence relation in terms of generating
functions with the use of the initial condition $C_{4}(v,u)=2uv^{2}$
(this holds immediately from the definitions), we can write
$$
\begin{array}{ll}
C(v,u;x) & =2uv^{2}x^{4}+\frac{vux}{1-u}(C(v,u;x)-C(vu,1;x))+\frac{vx}{1-v}%
(C(1,vu;x)-C(v,u;x)) \\
& +\frac{vx}{1-v}(B(vu;x)-vux^{3})-\frac{v}{1-v}(B(u;vx)-v^{3}ux^{3}) \\
& +\frac{vux}{1-u}(A(v;x)-vx^{2}-v(1+v)x^{3})-\frac{vux}{1-u}%
(A(vu;x)-vux^{2}-vu(1+vu)x^{3}) \\
& -\frac{2v^{4}ux^{5}}{(1-2vx)(1-u)}+\frac{2v^{4}u^{4}x^{5}}{(1-2vux)(1-u)}.%
\end{array}%
$$%
It is well known that this type of functional equations with several
variables are in general very hard to solve (see \emph{e.g.}
\cite{Mi03}). However, in our case, we are able to find an explicit
formula for the ordinary generating function $C(1,1;x)$, as it is
described below.

\textbf{Explicit formula for $C(1,1;x)$}: Substituting $u=v^{-1}$ in
the above functional equation gives
$$
\begin{array}{l}
C(v,v^{-1};x) =2vx^{4}-\frac{vx}{1-v}(C(v,v^{-1};x)-C(1,1;x))+\frac{vx}{1-v%
}(C(1,1;x)-C(v,v^{-1};x)) \\
 +\frac{vx}{1-v}(B(1;x)-x^{3})-\frac{x}{1-v}(B(v^{-1};vx)-v^{2}x^{3})-\frac{%
vx}{1-v}(A(v;x)-vx^{2}-v(1+v)x^{3}) \\
 +\frac{vx}{1-v}(A(1;x)-x^{2}-2x^{3})+\frac{2v^{4}x^{5}}{(1-2vx)(1-v)}-%
\frac{2vx^{5}}{(1-2x)(1-v)}.%
\end{array}%
$$%
This is equivalent to
$$
\begin{array}{l}
\left( 1+\frac{2vx}{1-v}\right)
C(v,v^{-1};x)=-(1+x+vx)vx^{3}+\frac{2vx}{
1-v}C(1,1;x)+\frac{vx}{1-v}B(1;x)-\frac{x}{1-v}B(v^{-1};vx) \\
\qquad-\frac{vx}{1-v}A(v;x)+\frac{vx}{1-v}A(1;x)+\frac{2v^{4}x^{5}}{(1-2vx)(1-v)}-\frac{2vx^{5}}{(1-2x)(1-v)},
\end{array}%
$$%
By taking $v=\frac{1}{1-2x}$ and using (\ref{eqff5}) and
(\ref{eqff6}),
$$
C(1,1;x)=\frac{2(3x-1)x^{2}}{(1-4x)^{3/2}}+\frac{2x^{2}(1-7x+15x^{2}-8x^{3})%
}{(1-2x)/(1-4x)^{2}}.
$$

\textbf{Explicit formula for $Q(1,1;x)$}: Equations (\ref{eqff3}), (\ref%
{eqff5}) and (\ref{eqff6}) give an explicit formula for $Q^{\prime }(1,1;x)$%
, namely $$Q^{\prime }(1,1;x)=\frac{(1-3x)x^{2}}{(1-4x)^{2}}-\frac{2x^{3}}{%
(1-4x)^{3/2}}.$$ Hence, by (\ref{eqff2}), we obtain that
$Q(1,1;x)=2Q^{\prime }(1,1;x)$ and the ordinary generating function
for the number of square permutations of length $n$ is given by
$1+x+2Q^{\prime }(1,1;x)$ ($1$ for the empty permutation and $x$ for
the permutation of length $1$), as required.
\end{proof}

\begin{corollary}
The number of polygons on $n$ vertices with four faces such that the
sources of the polygon lies on exactly two faces is given by
$$
2(n+2)4^{n-3}-2(n+1){{2n-4}\choose{n-3}}.
$$
\end{corollary}

\begin{proof}
The formula is obtained directly from Theorem \ref{hhs4} and Corollary \ref%
{co33}.
\end{proof}

\section{Open problems}

In this paper we have used a technique based on the kernel method to solve
functional equations for enumerating $k$-faces polygons on $n$ vertices,
where $k=2,3,4$. The results suggest the following problems:

\begin{itemize}
\item The most important question in our context is to find an explicit
formula for the number of $k$-faces polygons on $n$ vertices for any $k$.

\item We found that the number of $2$-faces polygons on $n$ vertices is
given by $\frac{2}{n-1}{\binom{{2n-4}}{{n-2}}}-2$ and that the number of
polygons on $n$ vertices with at most $2$ faces equals $\frac{2}{n-1}{\binom{%
{2n-4}}{{n-2}}}$, which is twice of the $n-2$-th Catalan number (see \cite[%
A000108]{JS}). This result can be explained combinatorially by considering
the number of permutations of length $n$ that can be put in increasing order
on two parallel queues (see Exercise 6.19 (jj) in \cite{st}). Also, we
proved analytically that the number of triangular permutations of length $%
n+2 $ is given by ${\binom{{2n}}{{n}}}$. This numbers appears frequently in
mathematics (see \cite[A000984]{JS}). One question is to find bijections
between consecutive-minima polygons and other mathematical objects.

\item Can we find a combinatorial interpretations for the formula $%
2(n+2)4^{n-3}-2(n+1){\binom{{2n-4}}{{n-3}}}$, the number of square
permutations of length $n$.

\item Theorem \ref{ppm4} gives the generating function $p_{m}(1;x)$ for the
number of $m$-parallel permutations of length $n$ when $m\leq 5$. Can we
find an explicit formula for any $m$. The same can be asked for the
generating functions for the number of $m$-triangular permutations of length 
$n$ (see Theorem \ref{ggm4}).

\item All questions about geometric properties of consecutive-minima
polygons remain open (for example, maximal perimeter, maximal area, number
of different polygons up to symmetries, \emph{etc.}).
\end{itemize}

\bigskip

\textbf{Acknowledgments.} Part of this work has been carried out while S.
Severini was visiting the Institut H. Poincar\'{e} (IHP), during the
programme \textquotedblleft Quantum Information, Computation and Complexity
(January 2006 - April 2006)\textquotedblright .

%==================================================================


\begin{thebibliography}{99}
\bibitem{AA} M. H. Albert and M. D. Atkinson, The enumeration of simple
permutations, \emph{Journal of Integer Sequences} \textbf{6} (2003), Article
03.4.4.

\bibitem{B02} C. Banderier, M. Bousquet-M{\'{e}}lou, A.~Denise, P.~Flajolet,
D.~Gardy, and D.~Gouyou-Beauchamps, Generating functions for generating
trees, \emph{Discr. Math.} \textbf{246:1-3}, 2002, 29--55.

\bibitem{Mi03} M. Bousquet-M\'{e}lou, Four classes of pattern-avoiding
permutations under one roof: generating trees with two labels, \emph{Elect.
J. Combin.} \textbf{9:2}, 2003, \#R19.

\bibitem{GM} {F. Ghassan and T. Mansour}, Three letters pattern avoiding
permutations, scanning elements method, and functional equations, preprint.

\bibitem{go} M. C. Golumbic, \emph{Algorithmic graph theory and perfect
graphs}. Second edition. With a foreword by Claude Berge. Annals of Discrete
Mathematics, 57. Elsevier Science B.V., Amsterdam, 2004.

\bibitem{co} J. E. Goodman, J. Pach and E. Welzl (Eds.), \emph{Combinatorial
and computational geometry}. Mathematical Sciences Research Institute
Publications, 52. Cambridge University Press, Cambridge, 2005.

\bibitem{M1} T. Mansour, The enumeration of permutations whose posets have a
maximal element, \emph{Advances in Applied Mathematics}, to appear (2006).

\bibitem{M2} T. Mansour, On an open problem of Green and Losonczy: exact
enumeration of freely braided permutations, \emph{Discrete Mathematics and
Theoretical Computer Science} \textbf{6:2} (2004), 461--470.

\bibitem{Rabbit} T. Mansour, Combinatorial methods and recurrence relations
with two indices,\emph{Journal of Difference Equations and Applications},
2006.

\bibitem{JS} N. J. A. Sloane and S. Plouffe, The Encyclopedia of Integer
Sequences, Academic Press, New York (1995).

\bibitem{st} R. P. Stanley, \emph{Enumerative combinatorics}. Volumes 1-2.
With a foreword by Gian-Carlo Rota. Corrected reprint of the 1986 original.
Cambridge Studies in Advanced Mathematics, 49. Cambridge University Press,
Cambridge, 1997-1999.

\bibitem{Va} V. Vatter, Enumeration schemes,
http://math.rutgers.edu/~vatter/publications/wilfplus, preprint.

\bibitem{We95} J. West, Generating trees and the {C}atalan and {S}chr\"{o}%
der numbers, \emph{Discr. Math.} \textbf{146} (1995) 247--262.

\bibitem{Z} D. Zeilberger, A proof of Julian West's conjecture that the
number of two stack-sortable permutations of length $n$ is $%
2(3n)!/((n+1)!(2n+1)!)$, \emph{Discr. Math.} \textbf{102} (1992) 85--93.

\bibitem{ZZ} D. Zeilberger. \newblock Enumeration schemes, and more
importantly, their automatic generation. 
\newblock {\em Annals of
Combin.} {\bf 2} (1998) 185--195.
\end{thebibliography}
\end{document}